\documentclass[12pt]{article}
\usepackage{amsmath,graphicx}

\newcommand{\BEAS}{\begin{eqnarray*}}
\newcommand{\EEAS}{\end{eqnarray*}}
\newcommand{\BEA}{\begin{eqnarray}}
\newcommand{\EEA}{\end{eqnarray}}
\newcommand{\BEQ}{\begin{equation}}
\newcommand{\EEQ}{\end{equation}}
\newcommand{\BIT}{\begin{itemize}}
\newcommand{\EIT}{\end{itemize}}
\newcommand{\BNUM}{\begin{enumerate}}
\newcommand{\ENUM}{\end{enumerate}}

\newcommand{\BA}{\begin{array}}
\newcommand{\EA}{\end{array}}

\newcommand{\ie}{{\it i.e.}}
\newcommand{\etc}{{\it etc.}}

\newcommand{\reals}{{\mbox{\bf R}}}

{\begin{quote}}{\end{quote}}

\makeatletter
\long\def\@makecaption#1#2{
   \vskip 9pt
   \begin{small}
   \setbox\@tempboxa\hbox{{\bf #1:} #2}
   \ifdim \wd\@tempboxa > 5.5in
        \begin{center}
        \begin{minipage}[t]{5.5in}
        \addtolength{\baselineskip}{-0.95pt}
        {\bf #1:} #2 \par
        \addtolength{\baselineskip}{0.95pt}
        \end{minipage}
        \end{center}
   \else
    \hbox to\hsize{\hfil\box\@tempboxa\hfil}
   \fi
   \end{small}\par
}
\makeatother

\newcounter{oursection}

\newcounter{lecture}

\title{Log-Linear Dynamical Systems}
\author{Steven Diamond}

\begin{document}
\maketitle

\begin{abstract}
  We present \emph{log-linear dynamical systems},
  a dynamical system model for positive
  quantities. We explain the connection to linear dynamical systems
  and show how convex optimization can be used to identify and control
  log-linear dynamical systems.
  We illustrate system identification and control with an example
  from predator-prey dynamics.
  We conclude by discussing potential applications of the
  proposed model.
\end{abstract}

\section{Introduction}
% Define log-linear dynamical systems,
% show relation to linear dynamical systems.
\paragraph{Definition.}
Suppose we have a positive time series $x_1,x_2,\ldots \in \reals^n_{++}$.
The time series is governed by \emph{log-linear dynamics} if the following
holds:
\begin{equation}\label{update}
  (x_{t+1})_i = c_i (x_t)_1^{A_{i1}} \cdots (x_t)_n^{A_{in}}, \quad i=1,\ldots,n,
\end{equation}
where $A \in \reals^{n \times n}$ and $c \in \reals^n_{++}$.
In other words, each entry of $x_{t+1}$ is a monomial function of the entries of
$x_t$ \cite{boyd2007tutorial}.

Under the change of variables $(\hat{x}_t)_i = \log((x_t)_i)$, 
the time series is governed by the linear dynamical system \cite[Chapter~9]{boyd2018introduction}
\[
  \hat{x}_{t+1} = A\hat{x}_t + \hat{c},
\]
where $\hat{c}_i = \log(c_i)$.
It follows that any fixed point $\hat{x}_t$ of the log-linear dynamical system
will satisfy $\hat{x}_t = (I - A)^{-1}\hat{c}$.
% \hat{x}_{t+1} = \hat{x}_t
% \hat{x}_t = A\hat{x}_t + \hat{c}
% \hat{x}_t = (I - A)^{-1}\hat{c}

\paragraph{Noise.} We can introduce noise into log-linear dynamical systems
with the update
\begin{equation}\label{noise1}
  (x_{t+1})_i = c_i (x_t)_1^{A_{i1}} \cdots (x_t)_n^{A_{in}}(z_t)_i, \quad i=1,\ldots,n,
\end{equation}
where $z_t \in \reals^n_{++}$ is drawn from a multivariate log-normal
distribution (or any distribution such that $z_t$ is a positive random variable).
Let $(\hat{z}_t)_i = \log((z_t)_i)$.
For example, we could have
\begin{equation}\label{noise2}
  \hat{x}_{t+1} = A\hat{x}_t + \hat{c} + \hat{z}_t,
\end{equation}
where $\hat{z}_t \sim \mathcal{N}(0,\sigma^2 I)$.

\section{System identification and control}\label{sec:id-control}
% Explain control of log-linear dynamical systems.
A major advantage of linear dynamical systems over non-linear models
is that linear dynamical systems are easy to identify and control.
The same holds true for log-linear dynamical systems,
since they are related to linear dynamical systems by a change of variables.

\paragraph{System identification.} Suppose we have measurements
$x_1,\ldots,x_T$ from a log-linear dynamical system.
Then we can identify the dynamics
by solving the least-squares problem
\begin{equation}\label{sysid}
\begin{array}{ll}
\mbox{minimize}   & \sum_{t=1}^{T-1} \|\hat{x}_{t+1} - A\hat{x}_t -
                    \hat{c} \|^2_2
\end{array}
\end{equation}
for variables $A$ and $\hat{c}$.
The residuals $\hat{r}_t = \hat{x}_{t+1} - A\hat{x}_t - \hat{c}$
can be used to estimate the noise $\hat{z}_t$  \cite{ljung2001system}.

\paragraph{Control.}
A controllable log-linear dynamical system features
a control input sequence $u_1,u_2,\ldots \in \reals^m_{++}$
in addition to the sequence $x_1,x_2,\ldots\in \reals^n_{++}$.
Instead of following the update rule in equation~\eqref{update},
the system follows the update rule
\[
  (x_{t+1})_i = c_i (x_t)_1^{A_{i1}} \cdots (x_t)_n^{A_{in}}
  (u_t)_1^{B_{i1}}\cdots(u_t)_m^{B_{im}}, \quad i=1,\ldots,n,
\]
where $B \in \reals^{n \times m}$.
The linear dynamical system formulation is
\[
  \hat{x}_{t+1} = A\hat{x}_t + B\hat{u}_t + \hat{c},
\]
where $(\hat{u}_t)_i = \log\left( (u_t)_i \right)$.
(Of course we can add noise as in equations~\eqref{noise1} and \eqref{noise2}.)

We control the system by choosing the control inputs.
For example, we might choose $u_1,\ldots,u_T$ by
solving an optimization problem of the form
\begin{equation}\label{control}
\begin{array}{ll}
\mbox{minimize}   & \sum_{t=1}^{T} f_t(\hat{x}_{t+1}) + g_t(\hat{u}_t) \\
\mbox{subject to} & \hat{x}_{t+1} = A\hat{x}_t + B\hat{u}_t +
                    \hat{c}, \quad t=1,\ldots,T,
\end{array}
\end{equation}
where $f_1,\ldots,f_{T}$ and $g_1,\ldots,g_{T}$ are convex functions.
Problem~\eqref{control} can be solved efficiently via convex optimization \cite{BoV:04}.

\section{Example}
We apply the system identification approach described above
to a dataset tracking hare and lynx populations in Canada \cite{elton1942ten}.
The hare and lynx exhibit a predator-prey relationship,
which is traditionally modeled with the Lotka-Volterra equations \cite{murray}.
Modeling the predator-prey relationship as a log-linear dynamical system
is an interesting alternative approach,
which also captures the non-linear population dynamics.

We estimate $A$ and $\hat{c}$ by solving problem~\eqref{sysid}.
We obtain
\[
  A = \left[ \begin{array}{cc}
                 .74 & -.37 \\
                 .21 & .70
                 \end{array}
               \right], \quad
               c = \left [ \begin{array}{c}
                             2.0 \\
                             .23
                           \end{array}
                   \right ]
\]
The coefficients of $A$ show that the hare population $(x_t)_1$ grows
proportionally with the previous year's hare population and
inversely proportionally with the previous year's lynx population.
The lynx population $(x_t)_2$ grows proportionally with both
the previous year's hare and lynx population.

\begin{figure}
\begin{center}
\includegraphics[width=0.9\textwidth]{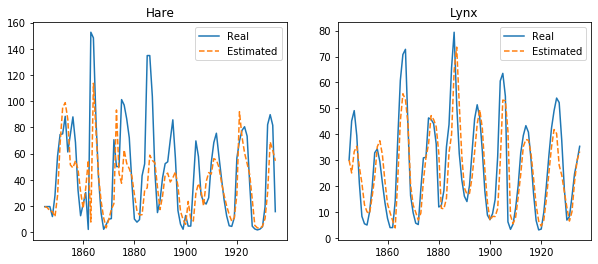}
\end{center}
\caption{Real and estimated hare and lynx populations (thousands of pelts taken).}
\label{estimates}
\end{figure}

The coefficients match our expectation that hares lead to more hares and
lynxes to fewer hares, while both hares and lynxes lead to more lynxes.
Figure~\ref{estimates} shows the true hare and lynx population (measured
indirectly via pelts taken) and the population predicted by equation~\eqref{update},
applied at each time step.
The fit is close enough to suggest the model is reasonable.

Identifying the log-linear dynamics of the predator-prey system
allows us to control the system
by introducing or removing hares and lynxes.
We could, for example, solve problem~\eqref{control} to drive the hare and lynx
populations along a desired trajectory,
or to keep the hare and lynx populations within desired bounds.
% show predator-prey example
% fit model, then control it

\section{Potential applications}
% chemical dynamics, any positive dynamical system like light intensity,
% population, energy, ...
Log-linear dynamical systems are a promising modeling paradigm.
Any positive dynamical system is a potential application.
Possibilities include human or animal populations, chemical concentrations
during reactions \cite{steinfeld1989chemical},
and light intensities on a sensor \cite{levoy2006light}.
We hope that researchers will find log-linear dynamical systems
a useful tool in their work.

\newpage
\appendix

\section{Modeling chemical reactions}

In this appendix we explore a possible application of log-linear dynamical
systems to chemical reactions.
Consider the rate equation for the reaction $A + B \rightleftharpoons AB$,
\ie, where species $A$ and $B$ combine to form $AB$:
\begin{equation}\label{eq:AB}
  \frac{\partial [AB]}{\partial t} = C[A][B],
\end{equation}
where $[A]$, $[B]$, and $[AB]$ are the molar concentrations of the reactants and
$C$ is a positive constant (in reality a function of temperature, \etc).

The rate equation looks similar to a log-linear dynamical system in that
chemical concentrations are always positive and the right-hand side is a
monomial function of the concentrations.
The main difference is that equation \eqref{eq:AB} is continuous,
a differential equation, while equation \eqref{update} is discrete.
We cannot arrive at equation \eqref{update} by discretizing equation
\eqref{eq:AB},
at least in any way obvious to the author.

Suppose though that we define the rate of change of $[AB]$ in a different way,
from the perspective of geometric rather than absolute growth.
Concretely, we consider the limit
\[
\lim_{\Delta \to 1} \left( \frac{[AB](\log t + \log \Delta))}{[AB](\log t)} \right)^{1/\log\Delta},
\]
which can be interpreted as the \emph{growth rate} of $[AB]$ at time $\log t$.
The limit is related to the derivative of $\log [AB]$.

We can redefine equation \eqref{eq:AB} in terms of the limit above, yielding
\begin{equation}
  \lim_{\Delta \to 1} \left( \frac{[AB](\log t + \log \Delta)}{[AB](\log t)} \right)^{1/\log\Delta} =
  C[A](\log t)[B](\log t). \label{chem-diff}
\end{equation}
We discretize equation \eqref{chem-diff} by fixing the value of $\Delta$.
We then have the log-linear dynamical system
\begin{equation}
  [AB](\log t + \log \Delta) =
  C^{\log\Delta}[A]^{\log\Delta}(\log t)[B]^{\log\Delta}(\log t)[AB](\log t). \label{ll-chem-diff}
\end{equation}

We have shown one way that chemical reactions could be modeled as a log-linear
dynamical system. Modeling chemical reactions as a log-linear dynamical system
would allow them to be easily controlled, as discussed in
Sec.~\ref{sec:id-control}.
The approach scales effortlessly to thousands or even millions of reactants.
Systems that control chemical reactions, such as cells,
could be modeled mathematically in a new and perhaps more sophisticated way.
% \[
%    \log[AB](\log t + \log \Delta) - \log[AB](\log t) = \log \Delta \left( \log C +
%     \log[A](\log t) + \log[B](\log t)  \right),
% \]

\newpage
\bibliography{llds}

\end{document}